\documentclass[12pt,a4paper]{article}
\usepackage{CJK}
\usepackage{indentfirst}
\usepackage[dvips]{graphics}
\usepackage{amsmath,amssymb,amsfonts,amsthm,graphics,CJK}
\usepackage{makeidx}
\begin{CJK}{GBK}{song}
\newtheorem{theorem}{Theorem}
\newtheorem{definition} {Definition}
\newtheorem{proposition}{Proposition}
\newtheorem{lemma}{Lemma}

\newtheorem{corollary}{Corollary}

\newtheorem{remark}{Remark}

\begin{document}
\title{\bf The (strong) rainbow connection numbers\\ of Cayley graphs of
Abelian groups\footnote{Supported by NSFC.}}
\author{\small Hengzhe Li, Xueliang Li, Sujuan Liu\\
\small Center for Combinatorics and LPMC-TJKLC\\
\small Nankai University, Tianjin 300071, China\\
\small lhz2010@mail.nankai.edu.cn; lxl@nankai.edu.cn;
sjliu0529@126.com}
\date{}
\maketitle

\begin{abstract}
A path in an edge-colored graph $G$, where adjacent edges may have
the same color, is called a rainbow path if no two edges of the path
are colored the same. The rainbow connection number $rc(G)$ of $G$
is the minimum integer $i$ for which there exists an
$i$-edge-coloring of $G$ such that every two distinct vertices of
$G$ are connected by a rainbow path. The strong rainbow connection
number $src(G)$ of $G$ is the minimum integer $i$ for which there
exists an $i$-edge-coloring of $G$ such that every two distinct
vertices $u$ and $v$ of $G$ are connected by a rainbow path of
length $d(u,v)$. In this paper, we give upper and lower bounds of
the (strong) rainbow connection Cayley graphs of Abelian groups.
Moreover, we determine the (strong) rainbow connection numbers of
some special cases.


\end{abstract}

\section{Introduction}

All graphs considered in this paper are undirected, finite and
simple. We refer to the book \cite{bondy} for graph theory notation
and terminology not described here. A path in an edge-colored graph
$G$, where adjacent edges may have the same color, is called a
$rainbow\ path$ if no two edges of the path are colored the same. An
edge-coloring of a graph $G$ is a $rainbow\ edge-coloring$ if every
two distinct vertices of $G$ are connected by a rainbow path.
Furthermore, a rainbow edge-coloring is a $strong\ rainbow\
edge-coloring$ if every two distinct vertices $x$ and $y$ of $G$ are
connected by a rainbow path with length $d(x,y)$. The $rainbow\
connection\ number\ rc(G)$ of $G$ is the minimum integer $i$ for
which there exists an $i$-edge-coloring of $G$ such that every two
distinct vertices of $G$ are connected by a rainbow path. If $G$ is
disconnected, we say that $rc(G)=0$ by convention. Furthermore, the
$strong\ rainbow\ connection\ number\ src(G)$ of $G$ is the minimum
integer $i$ for which there exists an $i$-edge-coloring of $G$ such
that every two distinct vertices $u$ and $v$ of $G$ are connected by
a rainbow path of length $d(u,v)$. It is easy to see that $D(G)\leq
rc(G)\leq src(G)$ for any connected graph $G$, where $D(G)$ is the
diameter of $G$.

The concept of rainbow connection number is of great use in
transferring information of high security in multicomputer networks.
Readers can see \cite{chartrand2} for details.

Let $\Gamma$ be a group, and let $a\in \Gamma$ be an element. We use
$\langle a\rangle$ to denote the cyclic subgroup of $\Gamma$
generated by $a$. The number of elements of $\langle a\rangle$ is
called the $order$ of $a$, denoted by $|a|$. A pair of elements $a$
and $b$ in a group commutes if $ab=ba$. A group is $Abelian$ if
every pair of its elements commutes.

\begin{definition}
The $Cayley\ graph$ of $\Gamma$ with respect to $S$ is the graph
$C(\Gamma,S)$ with vertex set $\Gamma$ in which two vertices $x$ and
$y$ are adjacent if and only if $xy^{-1}\in S$ (or equivalently,
$yx^{-1}\in S$), where $S\subseteq\Gamma\setminus{1}$ is closed
under taking inverse.
\end{definition}

It is well-known that $C(\Gamma,S)$ is connected if and only if $S$
is a generating set of $\Gamma$. The following conception is of
great convention in the sequel. An edge $xy$ is an $a-edge$ if
$xy^{-1}=a\in S$. It is not difficult to see that an $a$-edge $xy$
is also an $a^{-1}$-edge $yx$ since $S$ is closed under taking
inverse. Thus we do not distinguish $a$-edges and $a^{-1}$-edges in
the following arguments.

The hypercube is a well known model for computer networks which has
attracted many attentions in the past four decades, for example
\cite{bea,oh,park,saad,yang,yang2}. An $n$-$dimensional$ $hypercube$
is an undirected graph $Q_n=(V,E)$ with $|V|=2^n$  and
$|E|=n2^{n-1}$. Each vertex can be represented by an $n$-bit binary
string. There is an edge between two vertices whenever their binary
string representations differ in exactly one bit position.

\begin{definition}
The {\em recursive\ circulant} $G(N,d)$ has vertex set
$V=\{0,1,\ldots,N-1\}$, and edge set $E=\{vw\ |\ v\in V,w\in V |\
\mbox{ there exists}\ i,\ 0\leq i\leq \lceil\log_{d}N\rceil-1,\mbox{
such that}\ v+d^{i}\equiv$ $ w(\hspace{-5pt}\mod N)\}$.
\end{definition}

Combining Definition 1 and Definition 2, it can be seen that
$G(N,d)\cong Cay(Z_N,\{\pm d^{0},$ $\pm d^{1},\ldots,\pm
d^{\lceil\log_{d}N\rceil}\})$.

For the graph $G(rd^{m},d)$, we always assume $m\geq1$. Park and
Chwa showed in \cite{park} that the recursive circulant $G(N,d)$ has
a recursive structure when $N=rd^{m}, 1\leq c<d$ (See \cite{park}
for details), which is stated in the following theorem.

\begin{theorem}\cite{park}
Let $G(rd^{m},d)$ be a recursive circulant. Then
\begin{equation*}
D(G(rd^{m},d))=
\begin{cases}$$
\lfloor \frac {d}{2}\rfloor m + \lfloor r/2\rfloor & d\ is\ odd\\
\lfloor\frac{d-1}{2}m\rfloor+ \lfloor r/2\rfloor & both\ r\ and\ d\ is\ even\\
\lceil\frac{d-1}{2}m\rceil+ \lfloor r/2\rfloor & r\ is\ odd\ and\ d\
is \ even$$
\end{cases}
\end{equation*}
\end{theorem}

There have been some results on the (strong) rainbow connection
number of graphs, see \cite{chartrand,kri} for examples.

A $minimal\ generating\ set$ of a group $\Gamma$ is a generating set
$X$ such that no proper subset of $X$ is a generating set of
$\Gamma$. An $inverse\ closed\ minimal \ generating\ set$ of a group
$\Gamma$ is a set $X\cup X^{-1}$, such that $X$ is a minimal
generating set of $\Gamma$ and $X^{-1}=\{x\ |\ x\in X\}$. Clearly,
an inverse closed minimal generating set of $\Gamma$ contains only
one minimal generating set of $\Gamma$ if without distinguishing an
element $a$ and its inverse element $a^{-1}$.

In Section $2$, we show that $rc(C(\Gamma,S))\leq min\{\Sigma_{a\in
S^*} \lceil | a |/2\rceil\ |\ S^*\subseteq S$ is a minimal
generating set of $\Gamma$\}, where $\Gamma$ is an Abelian group.
Moreover, if $S$ is an inverse closed minimal generating set of
$\Gamma$, then $\Sigma_{a\in S^*} \lfloor | a |/2\rfloor\leq
rc(C(\Gamma,S))\leq src(C(\Gamma,S))\leq\Sigma_{a\in S^*} \lceil
|a|/2\rceil$, where $S^{*}\subseteq S$ is a minimal generating set
of $\Gamma$, and furthermore, if every element $a\in S$ has an even
order, then $rc(C(\Gamma,S))= src(C(\Gamma,S))=\Sigma_{a\in S^*}
|a|/2$.

In Section $3$, we show that $rc(G(rd^m,d))\leq src(G(rd^m,d))\leq
sm+\lceil r/2\rceil$, where $s$ is a constant related to $r$ and
$d$. Moreover, we prove that
$rc(G(r3^{m},3))=src(G(r3^{m},3))=m+\lfloor r/2 \rfloor$, and
$rc(G(2^{m},2))= src(G(2^{m},2))= m/2$ if $m$ is even.

\section{Cayley graphs of Abelian groups}

In this section, we first present the following elementary
proposition, then show the main result about Cayley graphs of
Abelian groups, and finally give some corollaries from the main
result.

\begin{proposition} If $H$ is a spanning subgraph of a graph $G$,
then $rc(G)\leq rc(H)$
\end{proposition}
It is easy to prove the above proposition since a rainbow
edge-coloring of $H$ induces a rainbow edge-coloring of $G$.

\begin{theorem}
Given an Abelian group $\Gamma$ and an inverse closed set
$S\subseteq \Gamma\setminus \{1\}$, we have the following results:

{\em  $(i)$} \[rc(C(\Gamma,S))\leq min\{\sum_{a\in S^*} \lceil | a
|/2\rceil\ |\ S^*\subseteq S\ is\ a\ \ minimal\ generating\ set\ of\
\Gamma\}.\]

{\em  $(ii)$} If $S$ is an inverse closed minimal generating set of
$\Gamma$, then\\
\[\sum_{a\in S^*} \lfloor | a |/2\rfloor\leq rc(C(\Gamma,S))\leq
src(C(\Gamma,S))\leq\sum_{a\in S^*} \lceil | a |/2\rceil,\]

\noindent where $S^{*}\subseteq S$ is a minimal generating set of
$\Gamma$.

Moreover, if every element $a\in S$ has an even order, then

\[
rc(C(\Gamma,S))= src(C(\Gamma,S))=\sum_{a\in S^*} |a|/2.\]

\end{theorem}
\begin{proof}
{\em  $(i)$} Note that a Cayley graph $C(\Gamma,S)$ is connected if
and only if $S$ is a generating set of $\Gamma$. Thus
$rc(C(\Gamma,S))=0$ if and only if $S$ is not a generating set of
$\Gamma$. Therefore, $(i)$ holds when $S$ is not a generating set of
$\Gamma$. Suppose $S$ is a generating set. We set
$\Gamma=\{v_1,v_2,\cdots,v_n\}$, and take any minimal generating set
$S^*=\{a_1,a_2,\cdots,a_r\}\subseteq S$ of $\Gamma$. Then
$C(\Gamma,S^{**})$ is a connected spanning subgraph of
$C(\Gamma,S)$, where $S^{**}=S^{*}\cup (S^{*})^{-1}$. It suffices to
show that $rc(C(\Gamma,S^{**}))\leq\Sigma_{a\in S^*} \lceil
|a|/2\rceil$ by Proposition $2$. For every $1\leq i\leq r$, we use
$M_i$ to denote the edge set of the $a_i$-edges, Then $M_i, 1\leq
i\leq r$, form a partition of $E(C(\Gamma,S^{**}))$.

Set $| a_i |=b_i$. If $b_{i}=2$, clearly, $M_i$ is a perfect
matching. Then we assign $M_{i}$ the color $(i,1)$. If $b_{i}\geq
3$, then we first pick up the identity element $u_{i,1}=1$ of
$\Gamma$. The vertex sequence $(1,a_{i}1,a^{2}_{i}1,\cdots,
a^{b_{i}}_{i}1=1)$ is a cycle, denoted by $C_{i,1}$. We second pick
up the vertex $u_{i,2}\in \Gamma$ such that $u_{i,2}\not\in
V(C_{i,1})$. Then $(u_{i,2},a_{i}u_{i,2},a^{2}_{i}u_{i,2},\cdots,
a^{b_{i}}_{i}u_{i,2}=u_{i,2})$ is a cycle, denoted by $C_{i,2}$. We
can successively do in this way until no vertex is left. Then we
obtain $n/b_{i}$ cycles $C_{i,1}, C_{i,2},\cdots,C_{i,n/b_{i}}$.
Now, color the edge of
$C_{i,k}=(u_{i,k},a_{i}u_{i,k},a^{2}_{i}u_{i,k},\cdots,
a^{b_{i}}_{i}u_{i,k}=u_{i,k}),1\leq i\leq r,1\leq k\leq
n/b_{i},b_{i}\neq 2$ by distinguish the following cases:

{\em Case $1.$} $b_{i}\geq 3$ is even.

Assign edges $(a^{j}_{i}u_{i,k})(a^{j+1}_{i}u_{i,k}),0\leq j\leq
b_{i}/2-1$ by colors $(i,j+1)$, and edges
$(a^{b_{i}/2+j}_{i}u_{i,k})$ $ ( a^{b_{i}/2+j+1}_{i}u_{i,k}),0\leq
j\leq b_{i}/2-1$ by colors $(i,j+1)$.

{\em Case $2.$} $b_{i}\geq 3$ is odd.

Assign edges $(a^{j}_{i}u_{i,k})$ $(a^{j+1}_{i}u_{i,k}),0\leq j\leq
(b_{i}-3)/2$ by color $(i,j+1)$, and edges
$(a^{(b_{i}+1)/2+j}_{i}u_{i,k})$ $
(a^{(b_{i}+1)/2+j+1}_{i}u_{i,k}),0\leq j\leq (b_{i}-3)/2$ by color
$(i,j+1)$. The edge $(a^{(b_{i}-1)/2}_{i}u_{i,k})$ $
(a^{(b_{i}+1)/2}_{i}u_{i,k})$ can have a color $(i,(b_{i}+1)/2)$.

Note that $C_{i,k},1\leq k \leq n/b_{i}$ need $\lceil b_{i}/2\rceil$
colors. Thus, the number of colors that we have used equals
$\Sigma_{a\in S^*} \lceil |a|/2\rceil$.

Next we will show that the above edge-coloring is a rainbow
edge-coloring $C(\Gamma,S^{**})$, that is, there exists a rainbow
path connecting any two distinct vertices $x,y$ of
$C(\Gamma,S^{**})$. we can assume that
$x=a^{i_{1}}_{1}a^{i_{2}}_{2}\cdots
a^{i_{r}}_{r},y=a^{j_{1}}_{1}a^{j_{2}}_{2}\cdots a^{j_{r}}_{r}$.
Clearly, $0\leq j_{k}-i_{k}$ $\leq \lceil
b_{k}/2\rceil(\hspace{-5pt}\mod\ b_{k})$ or $0\leq i_{k}-j_{k}\leq
\lceil b_{k}/2\rceil(\hspace{-5pt}\mod\ b_{k})$ for any $1\leq k
\leq r$, where ``$-a$" is ``$+a^{-1}$", without loss of generality,
we assume that $0\leq j_{k}-i_{k}\leq \lceil
b_{k}/2\rceil(\hspace{-5pt}\mod\ b_{k})$ for any $1\leq k \leq r$.
Then path $p=(x=a^{i_{1}}_{1}a^{i_{2}}_{2}\cdots
a^{i_{r}}_{r},a^{i_{1}+ 1}_{1}a^{i_{2}}_{2}\cdots
a^{i_{r}}_{r},\cdots,a^{j_{1}}_{1}a^{i_{2}}_{2}\cdots
a^{i_{r}}_{r},\cdots,a^{j_{1}}_{1}a^{j_{2}}_{2}\cdots
a^{j_{r}}_{r}=y)$ is a rainbow path between $x$ and $y$. This
completes the proof of part $(i)$.

For{\em  $(ii)$}, suppose $S$ is an inverse closed minimal
generating set of $\Gamma$. Note that $\Gamma$ has only one minimal
generating set $S^{*}$ contained in $S$ if without distinguishing an
element $a$ and its inverse element $a^{-1}$ and $S=S^{**}=S^{*}\cup
(S^{*})^{-1}$. It suffices to show that $D(C(\Gamma,S))=\Sigma_{a\in
S^{*}} \lfloor \frac{| a |}{2}\rfloor$ and the above edge-coloring
is a strong rainbow-coloring.

We first show that $D(C(\Gamma,S))=\Sigma_{a\in S^{*}} \lfloor
\frac{| a |}{2}\rfloor$. It is well-known that Cayley graphs are
vertex-transitive, we only consider the distance from $1$ to any
other vertex $x$ of $C(\Gamma,S)$. Without loss of generality,
assume that $x=a^{i_{1}}_{1}a^{i_{2}}_{2}\cdots a^{i_{r}}_{r}$
satisfying $i_{k}\leq  \lfloor b_{k}/2\rfloor, 1\leq k \leq r$.
Otherwise, we replace $a^{i_{k}}_{k}$ by
$(a^{-1}_{k})^{b_{k}-i_{k}}$ since
$a^{i_{k}}_{k}=(a^{-1}_{k})^{b_{k}-i_{k}}$ and $a^{-1}_{k}\in S$.
$P=(1=a^{0}_{1}a^{0}_{2}\cdots a^{0}_{r},a^{1}_{1}a^{0}_{2}\cdots
a^{0}_{r},\cdots, a^{i_{1}}_{1}a^{0}_{2}\cdots
a^{0}_{r},\cdots,a^{i_{1}}_{1}a^{i_{2}}_{2}\cdots a^{i_{r}}_{r}=x)$
is a path from $1$ to $x$ with length $\Sigma_{1\leq k\leq r}
i_{k}$. Thus, $D(C(\Gamma,S))$ $\leq\Sigma_{a\in S^{*}}  \lfloor | a
|/2\rfloor$. On the other hand, for $1$ and $x=a^{\lfloor
b_1/2\rfloor}_{1}a^{\lfloor b_2/2\rfloor}_{2}\cdots a^{\lfloor
b_r/2\rfloor}_{r}$, if $d(1,x)< \Sigma_{a\in S^{*}} \lfloor | a
|/2\rfloor$, then, by the pigeonhole principle, there exists an
integer $i$ such that the number of $a_{i}$-edges of the shortest
path from $1$ to $x$ less than $\lfloor b_i/2 \rfloor$, which is
impossible since $S$ is an inverse closed minimal generating set of
$\Gamma$. Therefore, $D(C(\Gamma,S))=\Sigma_{a\in S^{*}} \lfloor |a
|/2\rfloor$. Thus $\Sigma_{a\in S^*} \lfloor |a|/2\rfloor\leq
rc(C(\Gamma,S))\leq \Sigma_{a\in S^*} \lceil |a|/2\rceil$.

Next, we only need to show that for any $x,y\in V(C(\Gamma,S))$,
there exists a rainbow path with length $d(x,y)$ between $x$ and
$y$. We also can assume that $x=a^{i_{1}}_{1}a^{i_{2}}_{2}\cdots
a^{i_{r}}_{r},y=a^{j_{1}}_{1}a^{j_{2}}_{2}\cdots a^{j_{r}}_{r}$
satisfying $0\leq j_{k}-i_{k}\leq \lfloor
b_{k}/2\rfloor(\hspace{-5pt}\mod\ b_{k}), 1\leq k \leq r$. By a
similar argument of the diameter $D(\Gamma,S)$, we conclude that
$d(x,y)=\Sigma_{1\leq k\leq r} (j_k-i_k)$. Moreover, the path
$(x=a^{i_{1}}_{1}a^{i_{2}}_{2}\cdots
a^{i_{r}}_{r},a^{i_{1}+1}_{1}a^{i_{2}}_{2}\cdots
a^{i_{r}}_{r},\cdots, a^{j_{1}}_{1}a^{i_{2}}_{2}\cdots
a^{i_{r}}_{r},\cdots, a^{j_{1}}_{1}a^{j_{2}}_{2}\cdots
a^{j_{r}}_{r}=y)$ is a rainbow path from $x$ to $y$ with length
$d(x,y)=\Sigma_{1\leq k\leq r}  (j_k-i_k)$.

Now suppose every element $a\in S$ has an even order. Then
$\Sigma_{a\in S^*} \lfloor | a |/2\rfloor=\Sigma_{a\in S^*} \lceil
|a|/2\rceil =\Sigma_{a\in S^*} |a|/2$. So $rc(C(\Gamma,S))=
src(C(\Gamma,S))=\Sigma_{a\in S^*} |a|/2$. This completes the proof
of the theorem.
\end{proof}

The $Cartesian\ product$ of two simple graphs $G$ and $H$ is the
graph $G\Box H$ whose vertex set is $V(G)\times V(H)=\{(u,v)\ |\
u\in V(G),v\in V(H)\}$ and whose edge set is the set of all pairs
$(u_1, v_1)(u_2, v_2)$ such that either $u_1u_2\in E(G)$ and
$v_1=v_2$, or $v_1v_2\in E(H)$ and $u_1=u_2$.

Let $\mathbb{Z}_{n}$ be the $additive\ group$ of integers modulo
$n$. Note that $\mathbb{Z}^{n}_2$ is an Abelian group. The next
corollary follows from Theorem $2$ by taking
$\Gamma=\mathbb{Z}^{n}_2$ and $S=\{(1,0,\cdots,0)$,
$(0,1,\cdots,0),\cdots,(0,0,\cdots,1)\}$.

\begin{corollary}
Let $Q_{n}$ be an $n$-dimensional hypercube. Then
$$src(Q_{n})=rc(Q_{n})=n.$$
\end{corollary}

\begin{proof}
It is easy to see that $Q_{n}\cong P_{2}\Box P_{2}\Box P_{2}\cong
C(\mathbb{Z}^{n}_2,S)$ by the definitions of $n$-dimensional
hypercube, Cartesian product and Cayley graph
$C(\mathbb{Z}^{n}_2,S)$. Note that $S$ is a minimal generating set
of $\mathbb{Z}_{n}$, and also an inverse closed minimal generating
set of $\mathbb{Z}_{n}$. Thus, $src(Q_{n})=rc(Q_{n})=n$ follows from
Theorem $2$ and the fact that $|(\underbrace
{0,\cdots,1}_{k},\cdots,0)|=2$.
\end{proof}

$\mathcal{Z}=\mathbb{Z}_{n_1}\times
\mathbb{Z}_{n_2}\times\cdots\times \mathbb{Z}_{n_r}$ is an Abelian
group, where $n_k\geq 2,1\leq k\leq r$, and has an inverse closed
minimal generating set $S=\{(1,0,\cdots,0),(0,1,\cdots,0),
\cdots,(0,0,\cdots,1),
(n_1-1,0,\cdots,0),(0,n_2-1,\cdots,0),\cdots,(0,0,\cdots,n_r-1)\}$.
In the following corollary, by convenience, we set $C_{2}=P_{2}$.

\begin{corollary}
Let $C_{n_{k}}, n_{k}\geq 2,1\leq k\leq r$ be cycles, then

\[\sum_{1\leq k\leq r} \lfloor n_{k}/2\rfloor \leq rc(C_{n_1}\Box
C_{n_2}\Box C_{n_r})\leq src(C_{n_1}\Box C_{n_2}\Box
C_{n_r})\leq\sum_{1\leq k\leq r} \lceil n_{k}/2\rceil.\]

Moreover, if $n_k$ is even for every $1\leq k\leq r$, then

\[rc(C_{n_1}\Box C_{n_2}\Box C_{n_r})=src(C_{n_1}\Box C_{n_2}\Box
C_{n_r})=\sum_{1\leq k\leq r}  n_{k}/2.\]
\end{corollary}
\begin{proof}
We have $C_{n_1}\Box C_{n_2}\Box C_{n_r}\cong C(\mathcal{Z},S)$ by
the definitions of Cartesian product and Cayley graph
$C(\mathcal{Z},S)$. Moreover, $\mathcal{Z}$ has only one minimal
generating set
$S=\{(1,0,\cdots,0),\\(0,1,\cdots,0),\cdots,(0,0,\cdots,1)\}$
contained in $S$ if without distinguishing an element $a$ and its
inverse element $a^{-1}$, and $|(\underbrace
{0,\cdots,1}_{k},\cdots,0)|=n_{k},1\leq k\leq r$. Thus the first
inequality holds by Theorem $2$. Furthermore, if $n_k$ is even for
every $1\leq k\leq r$, we immediately deduce $src(C_{n_1}\Box
C_{n_2}\Box C_{n_r})=rc(C_{n_1}\Box C_{n_2}\Box
C_{n_r})=\Sigma_{1\leq k\leq r} n_{k}/2$ by Theorem $2$.
\end{proof}

\section{Recursive circulants}

Note that $G(rd^{m},d)\cong C(\mathbb{Z}_{rd^m},S)$, where $S=\{\pm
d^{0},$ $\pm d^{1},\ldots,\pm d^{\lceil\log_{d}rd^m\rceil}\}$.
Because $\{1\}$ is the only inverse closed minimal generating set
contained in $S$ of $\mathbb{Z}_{cd^{m}}$, we have
$rc(G(rd^{m},d))\leq src(G(rd^{m},d))\leq \lceil rd^{m}/2\rceil$ by
Theorem $3$. However, $\lceil rd^{m}/2\rceil$ is much larger than
$D(G(rd^{m},d))$ by Theorem $1$. So it is necessary to investigate
it further.

In this section, we first present some useful notations and one
helpful lemma from $\cite{park}$, then give upper and lower bounds,
and finally determine the (strong) rainbow connection numbers of
some special cases.

For $v\in V(G(rd^{m},d))$, a path between vertices $0$ and $v$ is a
sequence of vertices $v_0=0,v_1,\ldots,v_t=v$. Alternatively, it
also can be expressed by another sequence $b_1,b_2,\ldots,b_t$,
where $b_i=v_j-v_{j-1},\ 1\leq i\leq t$. The $ith$ vertex $v_i$ is
$\Sigma_{1\leq i\leq t}a_j$. It is easy to see that $b_i$ is either
$+d^j$ or $-d^j$ for some $j$. The following lemma is of great use.
We distinguish $+d^j$-edges and $-d^j$-edges in the following three
lemmas.

\begin{lemma}\cite{park} Let $P=b_1,b_2,\ldots,b_t$ be a shortest path from
$0$ to $v$.

$(i)$ $P$ does not have a pair of elements $+d^j$ and $-d^j$ for any
$j$.

$(ii)$ $P$ has less than $d\ ``+d^j$'s", and also has less than $d\
``-d^j$'s", for any $j$.
\end{lemma}

In fact, we can improve the above Lemma by the following three
lemmas.

\begin{lemma} For any vertex $v\in V(G(rd^{m},d)$, there exists a shortest
path $P=b_1,b_2,\ldots,b_t$ from $0$ to $v$ satisfying:

$(i)$ $P$ does not have a pair of elements $+d^j$ and $-d^j$ for any
$j$.

$(ii)$ $P$ has less than $\lfloor d/2\rfloor\ ``+d^j$'s", and also
has less than $\lfloor d/2\rfloor ``-d^j$'s", for any $0\leq j\leq
m-1$.

$(iii)$ $P$ has less than $\lfloor r/2\rfloor\ ``+d^{m}$'s", and has
less than $\lfloor r/2\rfloor\ ``-d^{m}$'s".
\end{lemma}
\begin{proof}
By Lemma $1$, we know that $(i)$ holds and $P$ has less than $d\
``+d^j$'s", and also has less than $d\ ``-d^j$'s", for any $j$. We
now show that $(ii)$ and $(iii)$ hold. Let path $P$ contains $r_j\
``+d^j$'s", or $r_j\ ``-d^j$'s". If $r_j\leq \lfloor d/2\rfloor$ for
all $0\leq j\leq m-1$ and $r_{m}\leq \lfloor r/2\rfloor$, we are
done. Otherwise, let $i$ be the smallest integer such that $i$ does
not satisfy the requirements of the lemma. Furthermore, without loss
of generality, we assume that $P$ contains $r_i\ ``+d^i$'s". We
consider the following two cases, according to whether $r\neq1$ or
$r=1$.

{\em Case $1.$} $r\neq 1$.

If $i\neq m$, then $r_i``+d^i$'s" can be replaced by one
$``+d^{i+1}$'s" and $(d-r_i)``-d^i$'s". Therefore, we obtain another
path $P'$ with length not larger than the length of $P$, and the
smallest integer $i$ which does not satisfies the requirements of
the lemma becomes larger. We can go on in this way until $i=m$. If
$r_{m}\leq \lfloor r/2\rfloor$, we are done. Otherwise, without loss
of generality, we assume that $P$ contains $r_m\ ``+d^m$'s". Then
$r_{m}\ ``+d^{m}$'s" can be replaced by $(r-r_{m})\ ``-d^{m}$'s". We
obtain a path $P^*$ with length not larger than the length of $P$,
and $P^*$ satisfies the requirements of the lemma.

{\em Case $2.$} $r=1$.

This case is similar to Case $1$ except that $G(d^{m},d)$ has no
$\pm d^m$-edge.

By this all possibilities have been exhausted and the proof is thus
complete.
\end{proof}

\begin{lemma} For any vertex $v$ of
$G(2^{m},2)$, there exists a shortest path $P=b_1,b_2,\ldots,b_t$
from $0$ to $v$ satisfying:

$(i)$ $P$ does not have a pair of elements $+2^j$ and $-2^j$ for any
$j$.

$(ii)$ For any integer $0\leq i\leq m-2$, $P$ does not have $``\pm
2^i$'s" or $``\pm 2^{i+1}$'s".
\end{lemma}
\begin{proof}
$(i)$ holds by Lemma $1$. Now let $i$ be the smallest integer such
that $P$ has both $``\pm 2^i$'s" and $``\pm 2^{i+1}$'s", and $P$
does not have $``\pm 2^j$'s" or $``\pm 2^{j+1}$'s" for all $j<i$. We
consider the following four cases.

{\em Case $1.$} $P$ has $+2^i$ and $+2^{i+1}$.

Then, $+2^i$ and $+2^{i+1}$ can be replaced by $+2^{i+2}$ and
$-2^i$.

{\em Case $2.$} $P$ has $+2^i$ and $-2^{i+1}$.

Then, $+2^i$ and $-2^{i+1}$ can be replaced by $-2^i$.

{\em Case $3.$} $P$ has $-2^i$ and $+2^{i+1}$.

Then, $-2^i$ and $+2^{i+1}$ can be replaced by $+2^i$.

{\em Case $4.$} $P$ has $-2^i$ and $-2^{i+1}$.

Then, $-2^i$ and $-2^{i+1}$ can be replaced by $-2^{i+2}$ and
$-2^i$.

We can go on in this way until we construct a new path $P^*$ from
$0$ to $v$ satisfying $(i)$ and $(ii)$. Note that the length of
$P^*$ is not larger than that of $P$. The proof is thus complete.
\end{proof}
\begin{remark}
Note that $G(rd^m,d)$ is vertex-transitive, thus, for any two
distinct vertices $u$ and $v$ of $G(rd^m,d)$, the analogous results
also hold.
\end{remark}
\begin{theorem}
Let $G(rd^{m},d)$ be a recursive circulant.

$(i)$ If $c\geq 4$ and $d\geq 5$, then,
\begin{equation*}
rc(G(rd^{m},d))\leq src(G(rd^{m},d))\leq
\begin{cases}$$
\frac{d}{2}m + \lceil r/2\rceil & d\ is\ even\\
\lfloor \frac{d}{2}\rfloor m + \lceil r/2\rceil & d\ is\ odd\ and\ r=\lfloor \frac{d}{2}\rfloor\\
(\lfloor \frac{d}{2}\rfloor+r)m+\lceil r/2\rceil &  d\ is\ odd\ and\ r < \lfloor \frac{d}{2}\rfloor\\
rm+\lceil r/2\rceil & d\ is\ odd\ and\ r>\lfloor\frac{d}{2}\rfloor$$
\end{cases}
\end{equation*}

In particular, if $d=2r+1$ and $r$ is even, then,

\[rc(G(rd^{m},d))= src(G(rd^{m},d))=rm+r/2.\]

$(ii)$ If $r\leq 3$ and $d\geq 4$, then, \begin{equation*}
rc(G(rd^{m},d))\leq src(G(rd^{m},d))\leq
\begin{cases}$$
\frac{d}{2}m + \lfloor r/2\rfloor & d\ is\ even\\
\lfloor \frac{d}{2}\rfloor m + \lfloor r/2\rfloor & d\ is\ odd\ and\ r=\lfloor \frac{d}{2}\rfloor\\
(\lfloor \frac{d}{2}\rfloor+r)m+\lfloor r/2\rfloor &  d\ is\ odd\ and\ r < \lfloor \frac{d}{2}\rfloor\\
rm+\lfloor r/2\rfloor & d\ is\ odd\ and\
r>\lfloor\frac{d}{2}\rfloor$$
\end{cases}
\end{equation*}

In particular, if $r=3$ and $d=7$, then,

\[rc(G(rd^{m},d))= src(G(rd^{m},d))=3m+1.\]

$(iii)$ If $d=3$, then,

$$rc(G(r3^{m},3))=src(G(r3^{m},3))=m+\lfloor r/2 \rfloor.$$

$(iv)$ If $d=2$, then,

$$\lfloor m/2\rfloor\leq rc(G(2^{m},2))\leq src(G(2^{m},2))\leq\lceil
m/2 \rceil.$$

In particular, if $m$ is even, then,

$$rc(G(2^{m},2))=src(G(2^{m},2))=m/2.$$

\end{theorem}
\begin{proof}
$(i)$ Note that $G(rd^{m},d)\cong C(\mathbb{Z}_{rd^m},S)$ if $r\geq
4$, where $S=\{\pm d^{0},$ $\pm d^{1},\ldots,\pm d^{m}\}$, and the
order of $\pm d^{i}$ is $|\pm d^{i}|=rd^{m-i}$. We use $M_i, 1\leq
i\leq r$, to denote the edge sets of the $\pm d^i$-edges. Then $M_i,
1\leq i\leq r$, form a partition of $E(G(rd^m,d))$.

{\bf First step:} Color the edges of $M_m$.

$M_m$ is a set of cycles with length $r$ since $r\geq 4$. By means
of the similar method of Theorem $2$, we can color the $\pm
d^m$-edges by $\lceil r/2\rceil$ colors.

{\bf Second step:} Color the edges of $M_i, 0\leq i \leq m-1$.

First, pick up the identity element $u_{i,1}=1$ of $\Gamma$. The
vertex sequence $(u_{i,1},u_{i,1}+d^{i},u_{i,1}+2d^{i},\cdots,
u_{i,1}+rd^{m-i}d^{i}=u_{i,1})$ is a cycle, denoted by $C_{i,1}$.
Second, pick up the vertex $u_{i,2}\in \Gamma$ such that
$u_{i,2}\not\in V(C_{i,1})$. Then
$(u_{i,2},u_{i,2}+d^{i},u_{i,2}+2d^{i},\cdots,
u_{i,2}+rd^{m-i}d^{i}=u_{i,2})$ is another cycle, denoted by
$C_{i,2}$. We can go on in this way until no vertex is left. Then we
obtain $d^i$ cycles $C_{i,1}, C_{i,2},\cdots,C_{i,d^i}$. Thus
$M_i=\bigcup_{1\leq k\leq d^i}C_{i,k}$ and $C_{i,k}\cap
C_{i,k'}=\emptyset$ for $k\neq k'$. We color the edges of
$C_{i,k}=(u_{i,k},u_{i,k}+d^{i},u_{i,k}+2d^{i},\cdots,
u_{i,k}+rd^{m-i}d^{i}=u_{i,k}),1\leq k\leq d^i,0\leq i\leq m-1$ by
distinguish the following cases:

{\em Case $1.$} $d$ is even.

Set $rd^{m-i}=\frac{d}{2}t$. Assign the edges
$(u_{i,k}+\frac{ld}{2}d^i+jd^i)(u_{i,k}+\frac{ld}{2}d^i+(j+1)d^i),0\leq
l\leq t-1,0\leq j\leq d/2-1$ by colors $(i,j+1)$. Thus, in this
case, the number of colors that we used to color all edges of
$G(rd^m,d)$ is $\frac{d}{2}m + \lceil r/2\rceil$.

{\em Case $2.$} $d$ is odd.

Then $d=2\lfloor d/2 \rfloor +1$, By direct computing, we have

$$d^{i}\equiv 1 (\hspace{-8pt}\mod \lfloor d/2\rfloor).$$

Multiplying both sides by $r$, we know
\begin{equation*}
rd^i\equiv
\begin{cases}
0\ & r=\lfloor d/2 \rfloor\\
r\ (\hspace{-8pt}\mod\lfloor d/2 \rfloor)& r< \lfloor d/2 \rfloor\\
(r-\lfloor d/2\rfloor)\ (\hspace{-8pt}\mod\lfloor d/2 \rfloor)& r>
\lfloor d/2 \rfloor
\end{cases}
\end{equation*}

Set $rd^{m-i}=\frac{d}{2}t+r'$, where $r'=0$ when $r=\lfloor
d/2\rfloor$, or $r'=r$ when $r< \lfloor d/2\rfloor$, or
$r'=r-\lfloor d/2\rfloor$ when $r> \lfloor d/2\rfloor$. We
distinguish the following there subcases.

{\em Case $2.1.$} $r=\lfloor d/2\rfloor$.

Similar to the method of Case $1$, we can color the edges of $M_i$
by colors $(i,1),(i,2),\ldots,$ $ (i,\lfloor d/2\rfloor)$. Thus, in
this case, the number of colors that we used to color all edges of
$G(rd^m,d)$ is $\lfloor\frac{d}{2}\rfloor m + \lceil r/2\rceil$.

{\em Case $2.2.$} $r< \lfloor d/2\rfloor$.

Assign the edges $(u_{i,k}+\lfloor\frac{d}{2}\rfloor
ld^i+jd^i)(u_{i,k}+\lfloor\frac{d}{2}\rfloor ld^i+(j+1)d^i),0\leq
l\leq t-1,0\leq j\leq \lfloor\frac{d}{2}\rfloor-1$ by colors
$(i,j+1)$. Assign the edges $(u_{i,k}+\lfloor\frac{d}{2}\rfloor
td^i+jd^i)(u_{i,k}+\lfloor\frac{d}{2}\rfloor td^i+(j+1)d^i),,0\leq
j\leq r'-1$ by colors $(i,\lfloor\frac{d}{2}\rfloor+j+1)$. Thus, in
this case, the number of colors that we used to color all edges of
$G(rd^m,d)$ is $(\lfloor d/2\rfloor+r)m+\lceil r/2\rceil$.

{\em Case $2.3.$} $r> \lfloor d/2\rfloor$.

Similar to the method of Case $2.2$, we can color the edges of $M_i$
by colors $(i,1),(i,2),\ldots,$ $ (i,\lfloor d/2\rfloor+r')$. Thus,
in this case, the number of colors that we used to color all edges
of $G(rd^m,d)$ is $(\lfloor d/2\rfloor+r')m+\lceil
r/2\rceil=rm+\lceil r/2\rceil$.

Next we will show that the above edge-coloring is a strong rainbow
edge-coloring of $G(rd^m,d)$, that is, there exists a rainbow path
with length $d(x,y)$ connecting any two distinct vertices $x,y$ of
$G(rd^m,d)$. By Lemma $2$, the commutativity of Abelian groups and
the method of the above edge-coloring, we can see that there exists
such a path.

In particular, if $d=2r+1$ and $r$ is even, we have
$\lfloor\frac{d}{2}\rfloor m + \lfloor r/2\rfloor\leq
rc(G(rd^{m},d))\leq src(G(rd^{m},d))\leq\lfloor\frac{d}{2}\rfloor m
+ \lceil r/2\rceil$ by Lemma $2$. Thus $rc(G(rd^{m},d))=
src(G(rd^{m},d))=\lfloor\frac{d}{2}\rfloor m + r/2$ since $\lfloor
r/2\rfloor=\lceil r/2\rceil=r/2$ when $r$ is even. The proof of
$(i)$ is thus complete.

For $(ii)$, we change the edge-coloring as follows: First, we can
color the $d^{i}$-edges with the similar method of Case $(i)$ for
$0\leq i\leq m-1$. Next, if $r=1$, there exists no $\pm
d^{m}$-edges. If $r=2$, the $\pm d^{m}$-edges form a prefect
matching. If $r=3$, the $\pm d^{m}$-edges are a set of disjoint
$K_3$'s, that is, the complete graph of order $3$. Thus, we can
color the $\pm d^{m}$-edges by $\lfloor r/2\rfloor$ colors. By Lemma
$2$, the commutativity of Abelian groups and the method of the above
edge-coloring, we know that the above edge-coloring is a strong
rainbow edge-coloring.

For $(iii)$, we change the edge-coloring as follows: If $r=1$, color
the $3^i$-edges by colors $(i,1)$ for $0\leq i\leq m-1$. If $r=2$,
color the $3^i$-edges by colors $(i,1)$ for $0\leq i\leq m$. By
Lemma $2$, the commutativity of Abelian groups and the method of
above edge-coloring, this is a strong rainbow edge-coloring. Thus
$rc(G(r3^{m},3))=src(G(r3^{m},3))=m+\lfloor r/2 \rfloor$ by Theorem
$1$.

For $(iv)$, we change the edge-coloring as follows: Color the $\pm
2^{2i}$-edges and the $\pm 2^{2i+1}$-edges by colors $(i,1)$ for
$0\leq i\leq \lceil m/2\rceil$. By Lemma $3$, the commutativity of
Abelian groups and the method of the above edge-coloring, this is a
strong rainbow edge-coloring. Thus $\lfloor m/2\rfloor\leq
rc(G(2^{m},2))\leq src(G(2^{m},2))\leq\lceil m/2 \rceil$ by Theorem
$1$. Furthermore, if $m$ is even, we have $\lfloor m/2\rfloor=\lceil
m/2 \rceil$, thus $rc(G(2^{m},2))=src(G(2^{m},2))=m/2$.

The proof is thus complete.
\end{proof}

We will present the following remark to complete the paper.

\begin{remark} In Theorem $2$, we show that, given an Abelian
group $\Gamma$ and an inverse closed minimal generating set
$S\subseteq \Gamma\setminus \{1\}$ of $\Gamma$, if every element
$a\in S$ has an even order, then $$src(C(\Gamma,S))=rc(C(\Gamma,S))=
\sum_{a\in S^*} | a |/2,$$ \noindent where $S^{*}\subseteq S$ is a
minimal generating set of $\Gamma$.

What happens when some element $a\in S$ has an odd order?  That is
the following open problem: given an Abelian group $\Gamma$ and an
inverse closed minimal generating set $S\subseteq \Gamma\setminus
\{1\}$ of $\Gamma$, is it true that

$$src(C(\Gamma,S))= rc(C(\Gamma,S))=\sum_{a\in S^{*}} \lceil
|a|/2\rceil\ ?$$

\noindent where $S^{*}\subseteq S$ is a minimal generating set of
$\Gamma$.
\end{remark}

\end{CJK}
\end{document}